\documentclass[11pt]{article}
\usepackage[utf8]{inputenc}

\usepackage{amsmath}
\usepackage{amsfonts}
\usepackage{accents}
\usepackage{graphicx}
\usepackage{geometry}
\usepackage{cancel}
\usepackage{fancyhdr}
\usepackage{dsfont}
\usepackage{mathrsfs}
\usepackage{upgreek}
\usepackage{amsthm}
\usepackage{amssymb}
\usepackage{mathtools}
\usepackage{pdfpages}
\usepackage{titlesec}
\usepackage{hyperref}
\usepackage{nicematrix}
\usepackage[sort&compress,comma,numbers]{natbib}
\setcitestyle{square}
\usepackage{newtxtext,newtxmath}
\usepackage{chngcntr}

\renewenvironment{abstract}
 {\par\noindent\textbf{\abstractname.}\ \ignorespaces}
 {\par\medskip}
 \renewcommand\thesubsubsection{\arabic{subsubsection}.}
 \titleformat{\subsubsection}{\itshape}{\thesubsubsection}{1em}{\itshape}

\title{\textbf{Real-time aerodynamic load estimation for hypersonics \\ via strain-based inverse maps}}

\author{Julie V. Pham\footnote{Graduate Research Assistant, Department of Aerospace Engineering and Engineering Mechanics.}, \hspace{1pt} Omar Ghattas\footnote{Professor, Walker Department of Mechanical Engineering, Oden Institute for Computational Engineering and Sciences.}, \hspace{1pt} Noel T. Clemens\footnote{Professor, Department of Aerospace Engineering and Engineering Mechanics},  \hspace{1pt} and Karen E. Willcox\footnote{Professor, Department of Aerospace Engineering and Engineering Mechanics.}}
\date{\emph{The University of Texas at Austin, Austin, TX 78712, USA}}

\begin{document}

\maketitle

\begin{abstract}
This work develops an efficient real-time inverse formulation for inferring the aerodynamic surface pressures on a hypersonic vehicle from sparse measurements of the structural strain. The approach aims to provide real-time estimates of the aerodynamic loads acting on the vehicle for ground and flight testing, as well as guidance, navigation, and control applications. Specifically, the approach targets hypersonic flight conditions where direct measurement of the surface pressures is challenging due to the harsh aerothermal environment. For problems employing a linear elastic structural model, we show that the inference problem can be posed as a least-squares problem with a linear constraint arising from a finite element discretization of the governing elasticity partial differential equation. Due to the linearity of the problem, an explicit solution is given by the normal equations. Pre-computation of the resulting inverse map enables rapid evaluation of the surface pressure and corresponding integrated quantities, such as the force and moment coefficients. The inverse approach additionally allows for uncertainty quantification, providing insights for theoretical recoverability and robustness to sensor noise. Numerical studies demonstrate the estimator performance for reconstructing the surface pressure field, as well as the force and moment coefficients, for the Initial Concept 3.X (IC3X) conceptual hypersonic vehicle.  
\end{abstract}

\section{Introduction}

Modern aerospace vehicles demand robust flight control systems that can operate autonomously under dynamic, uncertain, and extreme environments. Central to this requirement is the ability to obtain accurate information about the aerodynamic state of the vehicle. This information enables the deployment of advanced predictive capabilities for autonomy, digital twins, and guidance, navigation, and control (GNC). These desired capabilities become more challenging under hypersonic flight conditions, where the harsh aerothermal environment is inhospitable to many sensing technologies. In particular, it is difficult to obtain information about aerodynamic loads, such as the surface pressure field on the vehicle. The high-temperature hypersonic environment often limits the use of external surface sensors. Other conventional sensing technologies include inertial measurement units, GPS, satellites, and optical lidar, but these may not provide direct or sufficiently accurate estimates of the aerodynamic quantities of interest, and often have external dependencies which may be denied in adversarial scenarios. In this work we develop a vehicle-as-a-sensor concept~\cite{PHAM2022}, where the deformation of the vehicle is used to infer the instantaneous aerodynamic surface pressures. Using this strategy to obtain accurate aerodynamic load estimates in real time can help inform GNC systems, enhance maneuverability, and improve reliability for hypersonics. 

Aerodynamic loads are characterized by their integrated quantities, namely, the force and moment coefficients. Estimation of these aerodynamic parameters is commonly achieved using filtering methods such as the extended Kalman filter and unscented Kalman filter \cite{CHOWDHARY2010106, garcia1997}. The measurement models rely on the availability of measurements of the air flow angles, such as the angle of attack and sideslip angle, which are difficult to measure accurately. To mitigate this, Morelli et al. reconstruct the air flow angles from inertial data \cite{Morelli2012} for estimating the aerodynamic parameters.  In this paper, we take the approach of using discrete strain measurements to infer the full surface pressure field, and consequently the aerodynamic force and moment coefficients. To do this, we must obtain a (inverse) mapping from strain to aerodynamic quantities of interest. The strain response induced by the aerodynamic loads is governed by the partial differential equations (PDE) of linear elasticity. This PDE underpins the strain measurement model for the inference task, leading to a PDE-constrained inverse problem. 

Solutions to PDE-constrained inverse problems have been formulated in both the deterministic \cite{ghattas_willcox_2021} and statistical (Bayesian) \cite{stuart2010inverse} settings, and have great utility in "outer-loop" problems such as aerospace design optimization \cite{martins2005, hicken2010, martins2013}. These solutions are challenging to compute for a number of reasons. First, the measurements from which we infer the solution are often noisy and sparse. Additionally, PDE solution operators for many engineering problems of interest are smoothing, meaning that the output of the physical process is smoother than the input. For example, the stress field in a problem governed by linear elasticity is smoother than the loading field that induced it. These characteristics result in an \emph{ill-posed} inverse problem. A further challenge is that these inverse solutions often require many evaluations of the high-fidelity physics simulations, which is computationally intractable for real-time sensing problems. However, the linearity of the elasticity PDE in our problem allows for the formulation of a statistical least-squares problem with linear constraints, which can be solved explicitly via the normal equations. The resulting estimator is a strain-based inverse map that can be pre-computed and queried rapidly in real-time. Notably, using the statistical formulation of the least-squares problem, the known uncertainty due to sensor noise in the observed data also leads to a quantifiable uncertainty for the estimated solution. Consequently, the estimator fully captures the elastic behavior of the hypersonic vehicle described by the high-fidelity physics, while achieving the desired goal of rapid evaluation and uncertainty quantification. 

The proposed inverse approach relies on the availability and amenability of discrete, sparse strain measurements to provide sufficient information about the structural state induced by the aerodynamic loads. Strain-based sensing, and especially fiber-optics, has been used for various applications including finite element load updating \cite{kapania2007}, aeroelastic shape recovery \cite{freydin2019}, as well as real-time control and monitoring applications \cite{Derkevorkian2013, pena2018, MAININI2017296}. Our work extends the application of strain-based sensing to recover the full aerodynamic pressure field, as well as aerodynamic force and moment parameters. Some additional modeling aspects of this work rely on computational fluid dynamics (CFD) simulations for hypersonics to produce prior information for the surface pressure field. We note that the advancement of both simulation capabilities \cite{candler2021, mcnamara2011, san2014, cesnik2011, klock2017} and measurement technologies \cite{clemens2021} have progressed greatly, and are key enabling factors for inverse problems which must leverage both resources.

An alternative approach to solving the inverse problem is to employ machine learning to train a rapid-to-evaluate inverse mapping from the measurements to the quantities of interest. The training data can come from high-fidelity physics simulations, which are queried over a large range of flight conditions, as well as experimental data. Our previous work \cite{PHAM2022} has employed physics simulations to train optimal classification trees (OCTs) \cite{bertsimas2017optimal} for the inverse mapping.  The OCTs demonstrated strong prediction accuracy, in addition to being interpretable and rapid to evaluate. Neural networks are another possible option \cite{meinicke2024} which are fast to evaluate and have powerful representative capacity for nonlinear mappings, but they often lack interpretability and may perform poorly in extrapolation. While the accuracy performance may be reasonable, these machine-learned inverse maps lack a strong theoretical foundation for quantifying uncertainty and ensuring robustness in the predictions. In this paper, we seek to offer the same evaluation speed as the machine learning tools while remaining physics-informed and providing reliable quantification of uncertainty. 

The remainder of this paper is outlined as follows: in Section~\ref{sec:methods}, we present the formulation of the proposed estimator and the corresponding analytical uncertainty quantification. Section~\ref{sec:application} presents the testbed problem for demonstrating this work, along with numerical results. Finally, we discuss the conclusions and future work in Section~\ref{sec:conclusion}.

\section{Real-time surface pressure estimation: Problem statement and approaches} \label{sec:methods}

In this section, we present the inverse problem statement and approaches. Section~\ref{sec:definition} describes the problem definition and the PDE forward model. Section~\ref{sec:formulation} formulates the solution to the inverse problem via a least-squares estimator, with data-driven modeling considerations outlined in Section~\ref{section:considerations}. Section~\ref{sec:analysis} provides a detailed analysis of the expected performance and uncertainty of the estimator in two possible settings.

\subsection{Problem Definition}\label{sec:definition}

\begin{figure}[h]
	\centering
	\includegraphics[width = 6.6in]{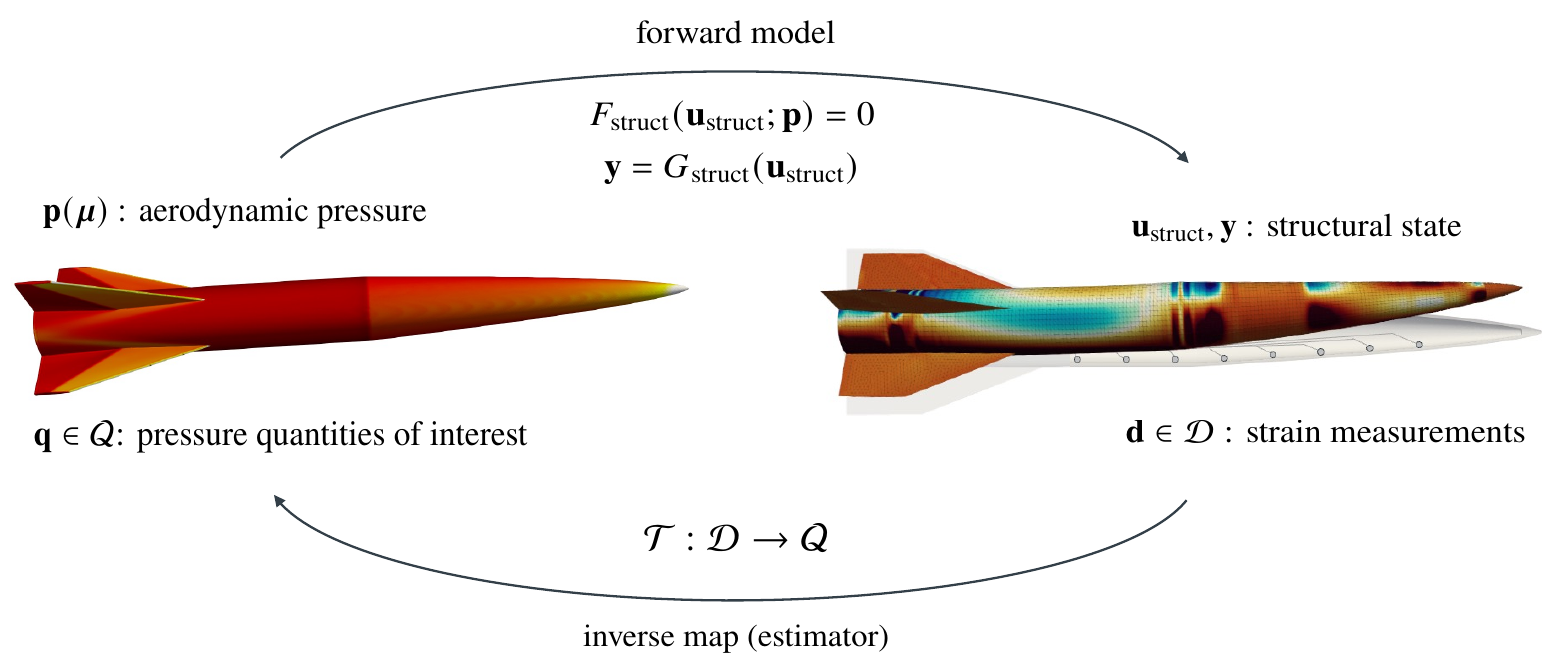}
	\caption{Inverse problem overview.}
	\label{fig:overview}
\end{figure}

We consider the structural and aerodynamic state of a hypersonic vehicle to be quasi-static at a particular instant in time. The structural state is observed through measurements of strain, with sensors placed on the interior surface of the vehicle. We denote the observed data by an $n_d$-dimensional vector $\mathbf{d} \in \mathcal{D}$, where $\mathcal{D} \subset \mathbb{R}^{n_d}$ is the space of possible strain measurements. The quantities of interest, $\mathbf{q} \in \mathcal{Q}$, are an $n_q$-dimensional discrete representation of the pressure field over the surface of the vehicle, where $\mathcal{Q} \subset \mathbb{R}^{n_q}$ is the space of possible pressure field representations. The goal is to infer $\mathbf{q}$ from the observed data $\mathbf{d}$. To perform the inference task, we seek an estimator that represents the inverse mapping 
\begin{equation}
    \mathcal{T} : \mathcal{D} \rightarrow \mathcal{Q}
\end{equation}
which can be evaluated rapidly in real-time. The relationship between the true aerodynamic surface pressure and the vehicle strain is described by a physics-based structural forward model, which we will later exploit in performing the inference task. Consider a discrete pressure field $\mathbf{p}(\boldsymbol{\mu})$ of dimension $n_p$ which is dependent on the flight conditions at a particular instant in time, $\boldsymbol{\mu} = [M, \alpha, \beta, H]^\top$, where $M$ is the Mach number, $\alpha$ the angle of attack, $\beta$ the sideslip angle, and $H$ the altitude. Our inverse problem formulation works with $\mathbf{q}$, which are the pressure quantities of interest corresponding to a discrete pressure field $\mathbf{p}(\boldsymbol{\mu})$. For example, $\mathbf{q}$ could be a vector containing $n_q < n_p$ nodal values of $\mathbf{p}$, or $\mathbf{q}$ could be a vector containing $n_q$ modal coefficients that define a representation of $\mathbf{p}$ in some modal basis. The pressure is the input to the structural forward model, denoted by $F_{\text{struct}}(\mathbf{u}_{\text{struct}}; \mathbf{p}) = 0$, where $\mathbf{u}_{\text{struct}}$ is the displacement of the deformed vehicle. After solving the forward model for $\mathbf{u}_\text{struct}$, the model strain data at the sensor locations is computed through the observation equation $\mathbf{y} = G_{\text{struct}}(\mathbf{u}_{\text{struct}})$. Ultimately, the model strain response $\mathbf{y}$ can  be compared to the measured strain data $\mathbf{d}$, thus driving the formulation of the inverse problem. Figure~\ref{fig:overview} illustrates the overview of the relationship between $\mathbf{p}$, $\mathbf{d}$, and intermediate quantities.

To explicitly define $F_{\text{struct}}$ and $G_{\text{struct}}$, we consider the equations of linear elasticity on the domain $\Omega \subset \mathbb{R}^3$, with exterior and interior surfaces $\partial\Omega_\text{ext}$ and $\partial\Omega_\text{int}$, respectively, and aft end $\partial\Omega_\text{aft}$.  The governing partial differential equation (PDE) and boundary conditions for the displacement $\mathbf{u}(\mathbf{x})$ for $\mathbf{x} \in \Omega$ are given by

\begin{equation} \label{eq:elasticity}
\begin{cases}
    \hspace{3mm} -\boldsymbol{\nabla} \cdot \boldsymbol{\sigma}(\mathbf{u}) = \mathbf{0} &\text{   in } \Omega  \\
    \hspace{6mm}\boldsymbol{\sigma}(\mathbf{u}) \cdot \mathbf{n} = \mathbf{0}  &\text{   on } \partial\Omega_\text{int}  \\
    \hspace{6mm}\boldsymbol{\sigma}(\mathbf{u})\cdot \mathbf{n} = \mathbf{t} &\text{   on } \partial\Omega_\text{ext} \\
    \hspace{12mm}\mathbf{u} = \mathbf{0} &\text{   on } \partial\Omega_\text{aft} \\
\end{cases}
\end{equation}
where $\boldsymbol{\sigma}$ is the Cauchy stress tensor, $\mathbf{n}$ is the outward-pointing surface normal, and $\mathbf{t}$ is the traction defined on the external boundary $\partial\Omega_\text{ext}$. Here, the traction boundary condition on $\partial\Omega_\text{ext}$ is derived from the aerodynamic surface pressure applied to the vehicle, and the Dirichlet condition on $\partial\Omega_\text{aft}$ is included for well-posedness of the PDE. The strain-displacement equation is given by $\boldsymbol{\varepsilon} = \frac{1}{2}(\nabla \mathbf{u} + \nabla \mathbf{u}^\top)$ where $\boldsymbol{\varepsilon}$ is the infinitesimal strain tensor, and the constitutive law is given by $\boldsymbol{\sigma} = 2 S \boldsymbol{\varepsilon}  + \frac{S(E-2S)}{3S-E} \text{tr}(\boldsymbol{\varepsilon})\mathbf{I}$, where $S$ is the shear modulus and $E$ is the Young's modulus of the material. We consider a finite element spatial discretization of the above elasticity equations. Upon discretization, the forward structural model is given by $F_\text{struct} = \mathbf{A}\mathbf{u}_\text{struct} - \mathbf{f} = \mathbf{0}$, where $\mathbf{A}$ is the stiffness matrix, and $\mathbf{f}$ are the nodal forces  on $\partial\Omega_\text{ext}$, which are computed from the surface pressure $\mathbf{p}$. We express the nodal forces as $\mathbf{f} = \mathbf{Cq}$, where $\mathbf{C}$ relates the pressure quantities of interest $\mathbf{q}$ to the nodal forces $\mathbf{f}$. The resulting forward structural model takes the form of the linear system $\mathbf{A}\mathbf{u}_\text{struct} - \mathbf{Cq} = \mathbf{0}$, where $\mathbf{A}$ is large, sparse square matrix with dimension $n_s$, the number of degrees of freedom in the structural discretization. We also define $\mathbf{y} = G_{\text{struct}}(\mathbf{u}_{\text{struct}}) = \mathbf{B}\mathbf{u}_{\text{struct}}$, where $\mathbf{B}$ is a spatial observation operator that (1) maps the discrete displacements to discrete strains using the strain-displacement relation and (2) extracts the particular strains corresponding to the direction and location of each sensor from the forward model.

\subsection{Real-time Inverse Solutions via Statistical Least Squares} \label{sec:formulation}

We formulate the inference problem of inferring pressure quantities of interest $\mathbf{q}$ from measurements $\mathbf{d}$ in general as a regularized, weighted least-squares problem. The main objective is to minimize the \emph{data misfit} between a given strain measurement and the forward model strain predictions at the sensor locations. We add regularization and weights to provide statistical information about the measurements, via the sensor noise covariance $\mathbf{\Gamma}_\text{n}$, as well as the quantities of interest, in the form of a prior mean and covariance, $\bar{\mathbf{q}}$ and $\mathbf{\Gamma}_\text{pr}$, respectively. The selection of these properties are further detailed in Section~\ref{section:considerations}. Then, the least-squares objective is given by 

\begin{equation}\label{eq:objective}
    \min_\mathbf{q} \hspace{1mm} \underbrace{\frac{1}{2}|| \mathbf{B}\mathbf{u}_{\text{struct}} - \mathbf{d} ||_{\Gamma_\text{n}^{-1}}^2}_\text{data misfit} + \underbrace{\frac{\gamma}{2}|| \mathbf{q} - \bar{\mathbf{q}} || _{\Gamma_\text{pr}^{-1}}^2}_\text{regularization} \text{\hspace{0.2in}subject to } \mathbf{A} \mathbf{u}_\text{struct} - \mathbf{C}\mathbf{q} = \mathbf{0},
\end{equation}
where $\gamma$ is a regularization parameter. The notable challenge is that the constraint is derived from a high-dimensional finite element model, and the solution of Equation~\ref{eq:objective} typically necessitates many evaluations of this expensive forward model. However, due to the linearity of the constraint, the least-squares problem offers an explicit solution via the normal equations, given by 

\begin{equation}\label{eq:ls-solution}
    \hat{\mathbf{q}} = (\mathbf{C}^\top \mathbf{A}^{-\top}\mathbf{B}^\top \mathbf{\Gamma}_\text{n}^{-1} \mathbf{B} \mathbf{A}^{-1}\mathbf{C} + \gamma \mathbf{\Gamma}_\text{pr}^{-1})^{-1} (\mathbf{C}^\top \mathbf{A}^{-\top} \mathbf{B}^\top \mathbf{\Gamma}_\text{n}^{-1} \mathbf{d} + \gamma \mathbf{\Gamma}_\text{pr}^{-1} \bar{\mathbf{q}}).
\end{equation}
This is the general form of our estimator, or inverse map, obtained via a least-squares approach. We note that $\mathbf{u}_\text{struct} = \mathbf{A}^{-1}\mathbf{Cq}$; it follows that the parameter-to-observable map is given by $\mathbf{y} = \mathbf{B}\mathbf{u}_\text{struct} = \mathbf{B}\mathbf{A}^{-1}\mathbf{Cq}$. We let $\mathbf{Z} =\mathbf{B} \mathbf{A}^{-1}\mathbf{C}$ denote the Jacobian of the parameter-to-observable map. Substituting in Equation~\ref{eq:ls-solution} and expanding, we obtain 
\begin{align} \label{eq:p2o-ls-solution}
    \hat{\mathbf{q}} &= (\mathbf{Z}^\top \mathbf{\Gamma}_\text{n}^{-1} \mathbf{Z} + \gamma \mathbf{\Gamma}_\text{pr}^{-1})^{-1} (\mathbf{Z}^\top \mathbf{\Gamma}_\text{n}^{-1} \mathbf{d} + \gamma \mathbf{\Gamma}_\text{pr}^{-1} \bar{\mathbf{q}})\\ \label{eq:segmented-ls-solution}
    &= \underbrace{[(\mathbf{Z}^\top \mathbf{\Gamma}_\text{n}^{-1} \mathbf{Z} + \gamma \mathbf{\Gamma}_\text{pr}^{-1})^{-1} \mathbf{Z}^\top \mathbf{\Gamma}_\text{n}^{-1}]}_{\mathbf{\coloneqq T}} \hspace{0.5mm} \mathbf{d} \hspace{1mm} + \hspace{1mm} \underbrace{(\mathbf{Z}^\top \mathbf{\Gamma}_\text{n}^{-1} \mathbf{Z} + \gamma \mathbf{\Gamma}_\text{pr}^{-1})^{-1} \gamma \mathbf{\Gamma}_\text{pr}^{-1}\bar{\mathbf{q}}}_{\coloneqq \mathbf{k}}. 
\end{align}
We note that the matrix $\mathbf{Z}$ contains the system matrices $\mathbf{A, B}$ and $\mathbf{C}$, which are measurement-independent, and can therefore be constructed in advance to pre-compute $\mathbf{Z}$. We additionally have a known, fixed $\mathbf{\Gamma}_\text{n}$, $\bar{\mathbf{q}}$, and $\mathbf{\Gamma}_\text{pr}$. Therefore, in Equation~\ref{eq:segmented-ls-solution}, the matrix $\mathbf{T} = (\mathbf{Z}^\top \mathbf{\Gamma}_\text{n}^{-1} \mathbf{Z} + \gamma \mathbf{\Gamma}_\text{pr}^{-1})^{-1} \mathbf{Z}^\top \mathbf{\Gamma}_\text{n}^{-1}$ can be entirely pre-computed, resulting in a matrix of size $n_q \times n_d$. Similarly, the term $\mathbf{k} =(\mathbf{Z}^\top \mathbf{\Gamma}_\text{n}^{-1} \mathbf{Z} + \gamma \mathbf{\Gamma}_\text{pr}^{-1})^{-1} \gamma\mathbf{\Gamma}_\text{pr}^{-1} \bar{\mathbf{q}}$ can be pre-computed. Upon deployment of the estimator, for a new measurement $\mathbf{d}$, we can rapidly estimate the quantities of interest $\mathbf{q}$ by simply performing the matrix-vector product $\mathbf{T}\mathbf{d}$, plus a shifting term $\mathbf{k}$ due to the prior mean. Notably, this computation can realistically be performed in real-time, even for large $n_q$, since the sensor dimension $n_d$ is typically moderate. In this case, \emph{real-time} refers to the typical measurement frequency of onboard sensors, which is on the order of 100 Hz; relative to this time scale, the matrix-vector product can be achieved an order of magnitude faster.  We note that integrated quantities of interest, such as force and moment coefficients, can also easily be computed by reconstructing the pressure field from $\mathbf{q}$, and may be especially useful for real-time guidance, navigation, and control applications. Consequently, this estimator is a possible efficient real-time sensing strategy for hypersonic vehicles.

\subsection{Data-Driven Modeling Considerations} \label{section:considerations}

\subsubsection{Noise and prior distributions}
We characterize two important features for statistical modeling: (1) the measurement noise model, and (2) the prior distribution. For the noise model, we consider additive Gaussian noise such that $\mathbf{d} = \mathbf{Bu}_\text{struct} + \boldsymbol{\eta} = \mathbf{y} + \boldsymbol{\eta}$, where the random variable $\boldsymbol{\eta} \sim (\mathbf{0}, \mathbf{\Gamma}_\text{n})$ is the noise discrepancy between the model strain response and the measurement. Here, we assume independent, identically distributed sensor measurements, with zero mean and variance $\sigma^2$, such that $\mathbf{\Gamma}_\text{n} = \sigma^2 \mathbf{I}$. This variance can be approximated through characterization of the physical sensor noise levels. Substituting $\mathbf{\Gamma}_\text{n}$ in Equation~\ref{eq:ls-solution}, the noise covariance-weighted least-squares solution simplifies to the ordinary least-squares solution if there is no regularization. In practice, we also perform a whitening transformation to normalize the variance $\boldsymbol{\eta}' = \mathbf{L}^{-1}\boldsymbol{\eta}$, where $\mathbf{L}$ is computed through the Cholesky factorization $\mathbf{\Gamma}_\text{n} = \mathbf{LL}^\top$ for numerical stability. While this is trivial for the diagonal covariance matrix in this setting, the factorization applies for more general covariance matrices.

The prior distribution of the quantities of interest is modeled by a Gaussian distribution with mean $\bar{\mathbf{p}}$ and covariance $\mathbf{\Gamma}_\text{pr}$, which are constructed as a pre-processing step. This is achieved in a data-driven fashion by using a computational fluid dynamics (CFD) model, denoted by $F_\text{aero}(\mathbf{p}; \boldsymbol{\mu}) = 0$, to compute surface pressure snapshots over a range of flight conditions $\boldsymbol{\mu}_j$ for $j = 1, ..., N$, where $N$ is the total number of flight condition combinations of interest. The pressure snapshots are denoted by $\mathbf{p}_j$, each of which correspond to the flight conditions $\boldsymbol{\mu}_j$. We define $\bar{\mathbf{p}}$ to be the mean over the snapshots. Then, we construct the covariance of the prior distribution by computing the covariance over the snapshots, 

\begin{equation}\label{eq:prior_cov}
    \boldsymbol{\Gamma}_\text{pr} = \sum_{j=1}^N \frac{(\mathbf{p}_j - \bar{\mathbf{p}})(\mathbf{p}_j - \bar{\mathbf{p}})^\top}{N-1}.
\end{equation}
This covariance captures spatial correlations in the pressure field from the pressure snapshots, including features such as smoothness or discontinuities. This is important for recovery of these features in the estimated pressure field if they are not strictly informed by the measurements, which we further describe in Section~\ref{sec:case2}.

\subsubsection{Dimensionality reduction}

Since the surface pressure field $\mathbf{p}$ is high-dimensional, we may seek to represent the pressure field in a low-dimensional subspace, which is one option for establishing well-posedness of the inverse problem. Here, we use the proper orthogonal decomposition (POD). The pressure snapshots described above, $\mathbf{p}_j$, are collected in a snapshot matrix $\mathbf{P}$. We center the snapshot matrix by the column-wise mean, $\bar{\mathbf{p}}$, to obtain $\tilde{\mathbf{P}}$. We perform the singular value decomposition on the centered snapshot matrix, $\tilde{\mathbf{P}} = \mathbf{V\Sigma W}^\top$. We then choose to retain the first $r$ left singular vectors. The reconstructed surface pressure field is then given by
\begin{equation}\label{eq:pod}
    \mathbf{p} \approx \mathbf{\bar{p}} + \sum_{i=1}^r c_i \mathbf{v}_i
\end{equation}
where $c_i$ and $\mathbf{v}_i$ are the $i$-th POD coefficient and  left singular vector, respectively. This allows us to represent the high-dimensional pressure field using $r$ coefficients.

\subsection{Analysis and Uncertainty Quantification} \label{sec:analysis}

Up to this point, the specific representation of pressure quantities of interest $\mathbf{q}$  have not been detailed. In the following sections, we outline two cases for $\mathbf{q}$: first, the case where we have fewer parameters than measurements, $n_q \leq n_d$, such that the least squares problem is over-determined, and second, the case where there are more parameters than measurements, $n_q > n_d$, such that the inverse problem is ill-posed. We provide a discussion of uncertainty quantification for each case, as well as an analysis of recoverable surface pressure information for a given sensor configuration. 
    
\subsubsection{Case 1: $n_q \leq n_d$}\label{sec:case1}

Consider the case where the number of parameters $n_q$ is less than or equal to the number of measurements, $n_d$. In this setting, we let $\mathbf{q}$ be a low-dimensional parameterization of the surface pressure field $\mathbf{p}$ so that $n_q \leq n_d$.  Here, we employ the POD decomposition given in Equation~\ref{eq:pod}. First, let $\mathbf{C}_\text{map}$ define the mapping from the pressure field to the resulting nodal forces, $\mathbf{f} = \mathbf{C}_\text{map}\mathbf{p}$.  We then define $\mathbf{C} \coloneqq \mathbf{C}_\text{map}\mathbf{V}_r$, where the columns of $\mathbf{V}_r$ are the first $r$ left singular vectors of the snapshot matrix, and $q_i \coloneqq c_i$, the corresponding POD coefficients. We seek to estimate the coefficients $\mathbf{q}$ from the data to reconstruct the pressure field.  Here, we consider the weighted least squares problem  without regularization ($\gamma=0$), meaning we seek the unique solution to the over-determined system of equations $\mathbf{Zq} = \mathbf{d}$, where $\mathbf{Z}$ is full-rank. To analyze this problem, we consider an infinitesimal perturbation $\delta \mathbf{d}$ to the strain response $\mathbf{y} = \mathbf{Zq}$. We seek to characterize the corresponding error in the least squares solution $\delta \mathbf{q} = \hat{\mathbf{q}} - \mathbf{q}$ with respect to the perturbation $\delta \mathbf{d}$.  Standard sensitivity analysis \cite{trefethenbau} of the ordinary least-squares solution to small perturbations in the right-hand side yields a relative condition number given by 

\begin{equation}
    \mathcal{K} = \sup_{\delta \mathbf{d}} \bigg( \frac{\|\delta \mathbf{q}\| }{\| \mathbf{q}\|} \bigg/ \frac{\|\delta \mathbf{d}\|}{\| \mathbf{d}\|} \bigg) = \frac{\kappa(\mathbf{Z})}{\nu \cos{\theta}} 
\end{equation}
where $\kappa(\mathbf{Z})$ is the condition number of the parameter to observable map, $\nu= \frac{\|\mathbf{Z}\|\|\mathbf{q}\|}{\|\mathbf{Zq}\|}$, and $\cos{\theta} = \frac{\|\mathbf{Zq}\|}{\|\mathbf{d}\|} $. For $\mathcal{K} >> 1$, the problem is ill-conditioned and we can expect large relative errors in the solution, even though $\delta \mathbf{d}$ is infinitesimal; this error becomes even larger once we consider a larger perturbation $\boldsymbol{\eta}$. We can improve the conditioning of the problem through the number of sensors as well as sensor placement, since the condition number is directly affected by the information gained from the sensors in the problem. In many least squares problems, it is also useful to add regularization so that the errors are not magnified due to the matrix conditioning. For our problem in Case 1, when $n_q \leq n_d$, instead we drive $\kappa(\mathbf{Z})$ to be sufficiently small through the sensor configuration, so we do not require regularization ($\gamma = 0$). To further quantify the uncertainty in the estimate, the covariance matrix of the estimator can be computed, given by the inverse of the Hessian of the data misfit term in our least squares objective function, 
\begin{equation}\label{eq:sls-cov}
    \mathbf{H}^{-1}_\text{misfit} = (\mathbf{C}^\top \mathbf{A}^{-\top}\mathbf{B}^\top \mathbf{\Gamma}_\text{n}^{-1} \mathbf{B} \mathbf{A}^{-1} \mathbf{C})^{-1}  = (\mathbf{Z}^\top \mathbf{\Gamma}_\text{n}^{-1} \mathbf{Z})^{-1}
\end{equation}
This is easy to compute for low-dimensional $\mathbf{q}$, and it provides an explicit quantification of the uncertainty on the estimated parameters.

\subsubsection{Case 2: $n_q > n_d$}\label{sec:case2}

We now consider the case where the number of parameters $n_q$ is larger than the number of measurements $n_d$. Particularly, we seek to estimate the distributed loads without restrictions to a particular low-dimensional parameterization using POD or otherwise. Instead, we seek to directly infer the surface load parameter field. In this setting, we define $\mathbf{q} \coloneqq \mathbf{p}$, and $\mathbf{C} \coloneqq \mathbf{C}_\text{map}$. This formulation results in an under-determined system of linear equations, for which there are many solutions. Of these solutions, it is sometimes sufficient to choose the minimum-norm solution. This solution lies within the row space of $\mathbf{Z}$, which is the data-informed subspace of the parameter field, and is the best we can achieve without prior information. However, we can significantly improve upon this result by constructing a prior distribution, which enters through the regularization term in the least-squares problem.

We compute the covariance matrix as in Equation~\ref{eq:prior_cov}. We note that this data-driven covariance will be rank-deficient when the number of snapshots $N$ used to construct the covariance is less than $n_q$ (the rank of this matrix will be at most $N$). Thus, we do not have $\boldsymbol{\Gamma}_\text{pr}^{-1}$ in the explicit solution in Equation~\ref{eq:p2o-ls-solution}, where it appears in the inverse of the Hessian  $\mathbf{H}^{-1} = (\mathbf{Z}^\top \mathbf{\Gamma}_\text{n}^{-1} \mathbf{Z} + \gamma\mathbf{\Gamma}_\text{pr}^{-1})^{-1} $ of the least-squares objective in Equation~\ref{eq:objective}. However, we can factor the inverse Hessian and apply the Woodbury matrix identity\footnote{The Woodbury identity may suffer from ill-conditioning depending on the construction of the regularization. It is important to verify the numerical stability of the computation to avoid this potential pitfall. }, resulting in  
\begin{equation}
    \mathbf{H}^{-1} = (\mathbf{Z}^\top \mathbf{\Gamma}_\text{n}^{-1} \mathbf{Z} + \gamma\mathbf{\Gamma}_\text{pr}^{-1})^{-1} = (\mathbf{\Gamma}_\text{pr}\mathbf{Z}^\top \mathbf{\Gamma}_\text{n}^{-1} \mathbf{Z} + \gamma \mathbf{I})^{-1} \mathbf{\Gamma}_\text{pr} = \frac{1}{\gamma}(\mathbf{\Gamma}_\text{pr} - \mathbf{\Gamma}_\text{pr}\mathbf{Z}^\top (\gamma\mathbf{\Gamma}_\text{n} + \mathbf{Z} \mathbf{\Gamma}_\text{pr} \mathbf{Z}^\top)^{-1} \mathbf{Z} \mathbf{\Gamma}_\text{pr} ).
\end{equation}
This allows for the use of a low-rank prior covariance, as well as cheap computation of $\mathbf{H}^{-1}$. Analysis of the  Hessian provides insights for pressure field recoverability through the prior distribution and the data. As demonstrated in \cite{ISAAC2015}, we solve the eigenvalue problem for non-zero eigenvalues of the \emph{prior-preconditioned} Hessian of the data misfit, formulated as
\begin{equation}\label{eq:eigproblem}
    \boldsymbol{\Gamma}_\text{pr} \mathbf{H}_\text{misfit} \mathbf{b}_k = \lambda_k \mathbf{b}_k
\end{equation}
where $\lambda_k$ and $\mathbf{b}_k$ are the $k$-th largest eigenvalue and corresponding eigenvector, respectively, and $\mathbf{H}_\text{misfit} = \mathbf{Z}^\top \boldsymbol{\Gamma}_\text{n}^{-1} \mathbf{Z}$ as previously defined.  The eigenvectors (modes) of the resulting low-rank spectral decomposition of the prior-preconditioned Hessian define the possible solution space of the least-squares estimator. In a noise-free setting, these modes can be perfectly recovered via the prior distribution and the data. The orthogonal complement, which corresponds to higher frequency (spatially oscillatory) modes, is eliminated from the possible solution space. Since we expect our parameter field to be mostly smooth, this elimination results in a physically tractable solution. Additionally, standard spectral analysis \cite{ghattas_willcox_2021} shows that the errors due to noise in the Fourier coefficient corresponding to the $k$-th mode will scale with $\frac{1}{\lambda_k}$, which leads to amplified errors when $\lambda_k$ is small. In turn, the regularization parameter $\gamma$ must be chosen to shift the eigenvalues to reduce the noise amplification, which then also eliminates the contribution of the modes corresponding to eigenvalues $\lambda_k << \gamma$.  We tune this parameter based on the expected sensor noise levels in order to balance estimator accuracy and noise mitigation.

\section{Application: IC3X Testbed Problem} \label{sec:application}

\begin{figure}[h]
	\centering
	\includegraphics[width = 5.5in]{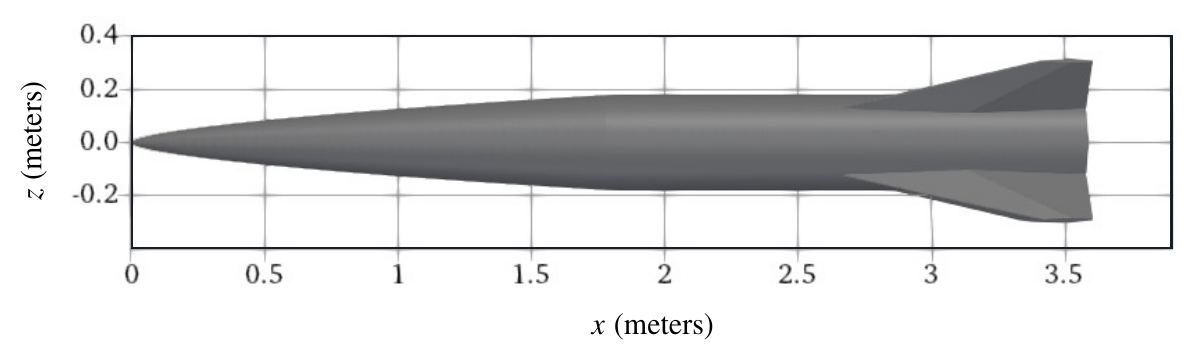}
	\caption{IC3X geometry.}
	\label{fig:ic3x_geom}
\end{figure}

To demonstrate the inverse approaches outlined above, we consider the Initial Concept 3.X (IC3X) hypersonic vehicle. The IC3X was initially proposed by Pasiliao et al. \cite{Pasiliao}, and a detailed finite element model for the vehicle was developed by Witeof et al. \cite{witeof2014initial}. The IC3X is a representative hypersonic vehicle with boost, cruise, and terminal phases in a nominal trajectory. The vehicle is propelled by a scramjet in the cruise phase. The geometry of the vehicle is shown in Figure~\ref{fig:ic3x_geom}. In this paper, we conduct an initial study using a modified version of the IC3X vehicle with removed fins (fuselage only), and a hollow internal structure for ease of mesh generation.  

\subsection{IC3X modeling}

 The forward aerodynamics for the IC3X are solved using CART3D~\cite{cart3d, cart, aftosmis1997cart3d}, which solves the compressible Euler equations for inviscid steady fluid dynamics. Given the vehicle surface geometry, CART3D employs a multilevel Cartesian cut-cell mesh with embedded boundaries, with adjoint-based adaptive mesh refinement \cite{nemec2007}. CART3D enables rapid databasing of surface pressure solutions over many different flight conditions with low computational cost, in comparison to higher fidelity physics models (i.e. Reynolds-Averaged Navier Stokes). We note the inviscid solver captures the surface pressure with high accuracy, but  does not capture the drag forces due to skin friction. Future work will consider estimation of the viscous drag component. 
 
 The structural model is constructed and solved using FEniCS\cite{AlnaesEtal2015, LoggEtal2012}, which is a software project providing a collection of packages, including DOLFIN\cite{LoggWells2010}, that can be used for solving PDEs with the finite element method. Specifically, we employ FEniCs for the explicit construction of the required system matrices $\mathbf{B}, \mathbf{A}, \mathbf{C}$ that appear in Equation~\ref{eq:ls-solution}.  We note that the estimator can also be implemented with access to only the system matrix actions on a vector, which may be needed for large-scale problems where it is computationally infeasible to store the system matrices explicitly. In this work, we do not consider thermal effects on the strain response of the vehicle due to thermal expansion, nor the temperature dependence of some material properties, which can potentially introduce nonlinearities. These considerations are to be addressed in future work.   
 
\begin{figure}[h]
	\centering
	\includegraphics[width = 5in]{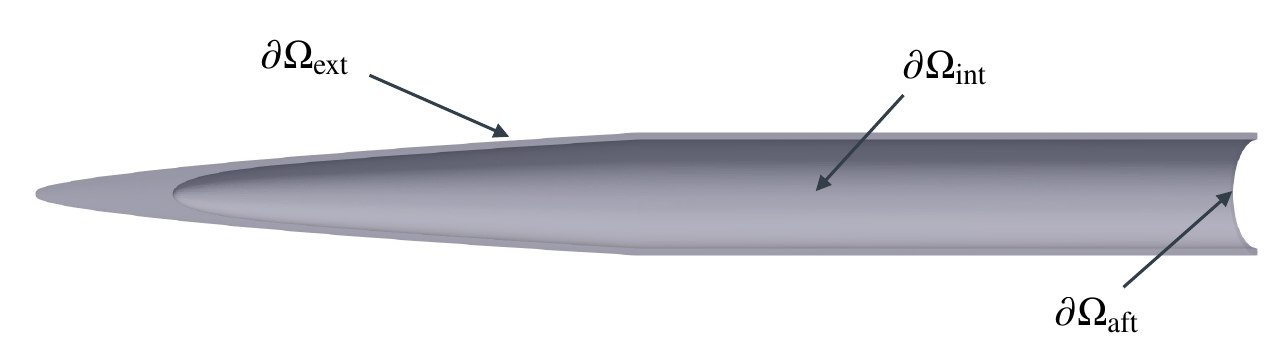}
	\caption{Boundary surfaces of the modified IC3X. }
	\label{fig:boundaries}
\end{figure}

 We produce the triangular surface mesh for CART3D and the tetrahedral volume mesh for the structural model using Gmsh \cite{gmsh}. The resulting surface mesh contains $n_p = 60{,}538$  degrees of freedom, and the resulting structural mesh contains $n_s =  461{,}664$ degrees of freedom. The boundary surfaces for the structural model, with boundary conditions given in Equation~\ref{eq:elasticity}, are illustrated in Figure~\ref{fig:boundaries}. 

\subsubsection{Sensor configuration and noise statistics}

\begin{figure}[h]
	\centering
	\includegraphics[width = 6.5in]{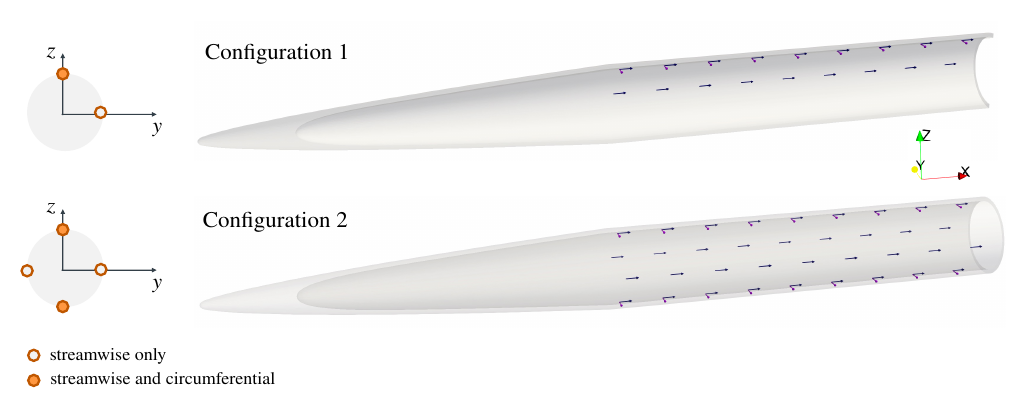}
        \caption{Strain sensor configurations. Left: yz-plane view indicating sensor row measurement direction. Right: 3D view of sensor locations, with arrows indicating measurement direction. }
	\label{fig:sensor-config}
\end{figure}

In this work, the strain sensor configurations consist of strain gauges placed in rows on the inner surface of the hollow IC3X vehicle structure. Each row consists of nine evenly spaced strain gauges. The in-plane strain of the inner surface can be measured in a particular direction. Measuring the strain in the stream-wise direction (i.e., x-direction) captures the bending strain, which dominates the strain response when the total angle of attack is non-zero. The circumferential strain can also be measured, capturing the strain due to hoop stresses. Figure~\ref{fig:sensor-config} depicts the sensor placement for two and four rows of sensors, which we will refer to as Configuration 1 and 2, respectively. Configuration 1 contains 27 sensors, of which 18 sensors measure the x-direction strain, and 9 sensors measure the circumferential strain. Configuration 2 contains 54 sensors, which are the same as Configuration 1 plus their mirrored locations through the x-axis (i.e., sensors in the xy-plane are reflected across the xz-plane, and sensors in the xz-plane are reflected across the xy-plane). The observation operator $\mathbf{B}$ is constructed for each of these configurations to map the forward model displacement to the strain in the corresponding direction and location as the physical sensors. As described in Section~\ref{section:considerations}, we assume independent, identically distributed Gaussian noise for each sensor with zero mean and standard deviation $\sigma$. For this study, we choose $\sigma$ to be equal to 1\% of the median strain magnitude over all sensors and flight conditions, based on estimated noise characteristics for physical strain gauges \cite{MAININI2017296, fiberoptic, meinicke_ifasd}. 

\subsubsection{Flight conditions and surface pressure}
\begin{figure}[h]
	\centering
	\includegraphics[width = 6.5in]{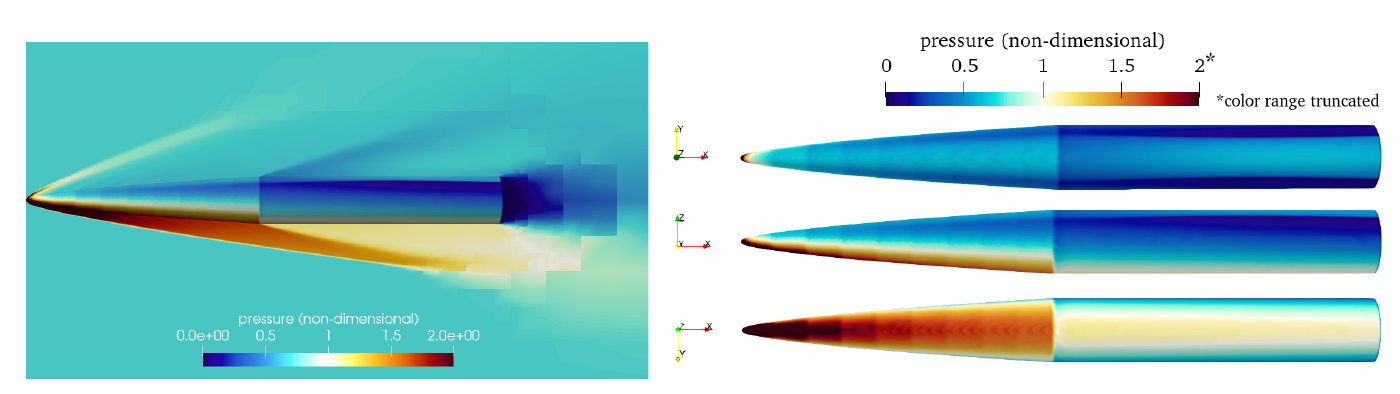}
	\caption{Pressure at conditions $M=5, \alpha=8, \beta=0$. Left: Side view of pressure over 2D slice of fluid domain. Right: (top) leeward, (middle) side, and (bottom) windward views of the surface pressure field.\protect\footnotemark} 
	\label{fig:pressure}
\end{figure}

For a nominal trajectory of the IC3X, we consider the following flight condition ranges of interest: Mach number 5 -- 7, angle of attack 0\textdegree -- 10\textdegree, sideslip angle 0\textdegree -- 10\textdegree, and altitude 20 km. For dimensionality reduction in Case 1, we solve the forward aerodynamics for the set $P_1$ of all combinations of the parameters $M = \{5.0, 5.5, 6.0, 6.5, 7.0\}$, $\alpha = \{0, 2, 4, 6, 8, 10\}$, $\beta = \{0, 5, 10\}$ at $H = 20$km altitude, and collect the corresponding pressure snapshots. This results in a total of $N = 90$ snapshots. Due to symmetry of the vehicle, we choose positive values of $\alpha$ and $\beta$ for computing the POD in Case 1 using the centered snapshots as in Equation~\ref{eq:pod}, which incorporates higher spatial frequency features in the first few modes. For the data-driven prior distribution in Case 2, we use the set $P_2$ of all combinations of the parameters $M=\{5, 6, 7\}$, $\alpha = \{-8, -6, -4, -2, 0, 2, 4, 6, 8\}$  and  $\beta = \{-8, -6, -4, -2, 0, 2, 4, 6, 8\}$, at $H = 20$km altitude. The corresponding $N=243$ snapshots are then used to compute the prior mean and covariance as in Equation~\ref{eq:prior_cov}. Figure~\ref{fig:pressure} gives a visualization of the pressure at flight conditions $M=5, \alpha=8$, and $\beta=0$, with a 2D slice of the fluid domain on the left, and the surface pressure on the right. We note that only the surface pressure is required, thus allowing a coarser mesh away from the surface as seen in the fluid domain slice. At the axial station $x=1.75$m where the angle of the fuselage changes, there is a clear discontinuity in the surface pressure, which is caused by the geometry and occurs at all flight conditions. Additionally, at larger angles of attack, there exists a discontinuity (shock) due to counter-rotating vortices that occur on the leeward side of the vehicle, leading to a region of higher surface pressure along the body. This is visible in the top right view of Figure~\ref{fig:pressure}. These types of features are common for hypersonic vehicles, and we seek to accurately estimate these surface pressure fields.  

\footnotetext{For visual clarity, the color range in this and subsequent pressure color bars is truncated at the indicated maximum. Higher pressures exist in a small region at the IC3X nose. }

\subsection{Case 1 : Numerical Results}

In this section, we provide numerical results for the scenario where $n_q \leq n_d$. The quantities of interest $\mathbf{q}$ we seek to estimate are the POD coefficients parameterizing the surface pressure. For the pressure snapshots corresponding to $P_1$, the cumulative energy retained by the first five modes is 99.1\%, thus we choose to retain $r=5$ singular vectors. Therefore, the parameter-to-observable map $\mathbf{Z}$ has dimensions $n_d \times r$. First, we demonstrate the impact of the sensor configuration on the conditioning of the least squares problem. Table~\ref{tab:relative-cond} shows the relative condition number $\mathcal{K}$ using $\|\mathbf{d}\| = 100\sqrt{n_d}$ and $\|\mathbf{q}\| = 1$ for each sensor configuration. Configuration 2 has the smaller condition number, since it has a larger number of sensors, and will therefore have less sensitivity to noise perturbations for estimating $\mathbf{q}$.

\begin{figure}[h!]
    \centering
    \includegraphics[width=5.2in]{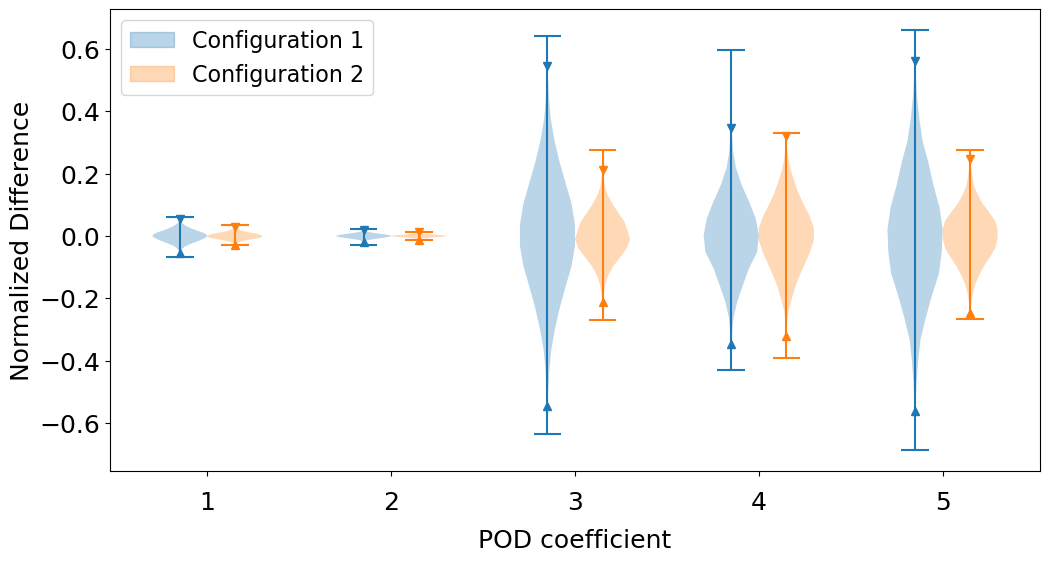}
    \caption{Normalized difference in estimated POD coefficients for Configuration 1 and 2 with $\mathbf{3}\boldsymbol{\sigma}$ uncertainty bounds denoted by the triangle markers.}
    \label{fig:pod_relative_errors}
\end{figure}

\begin{figure}[h!]
    \centering
    \includegraphics[width=3.9in]{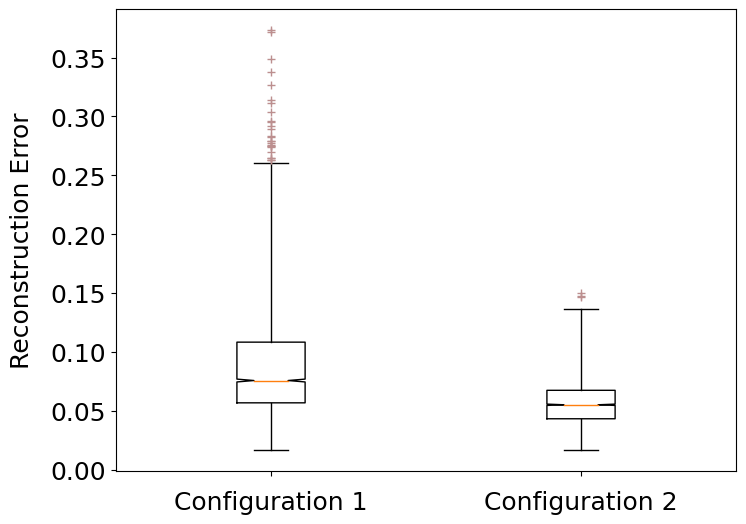}
    \caption{Relative L2 error in reconstructed pressure field for Configuration 1 and 2. }
    \label{fig:pod_recon_errors}
\end{figure}

\begin{table}[h]
\begin{center}
\caption{Relative condition number for sensor configurations.}
\label{tab:relative-cond}
\begin{tabular}{ c c c } 
 \hline\hline
  & Number of Sensors & Relative Condition Number $\mathcal{K}$  \\ 
 \hline
 Configuration 1 & 27 & 89.61 \\ 
 Configuration 2 & 54 & 47.96 \\ 
 \hline
\end{tabular}
\end{center}
\end{table}
To numerically demonstrate the performance of the estimator, we produce synthetic noisy measurements by sampling the noise model $\boldsymbol{\eta} \sim (\mathbf{0}, \mathbf{\Gamma}_n)$, and adding the realization of noise to the model strain response. The estimated POD coefficients $\hat{\mathbf{q}} = \mathbf{Td}$ with $\gamma=0$ are computed for each synthetic measurement. Each estimated coefficient $\hat{q}_i$ is compared to the reference POD coefficient $q_{i}$ at a particular flight condition using the normalized difference, computed by
\begin{equation}
    e^\text{POD}_i = \frac{\hat{q}_i - q_i}{\text{range}(q_i)}
\end{equation}
where $\text{range}(q_i)$ is the range of the $i$-th POD coefficient over all POD representations of the pressure snapshots in $P_1$. Figure~\ref{fig:pod_relative_errors} shows the normalized difference metric in the predicted coefficients over 50 samples from each flight condition in $P_1$, resulting in a total of 4500 samples. We also show the $3\sigma$ uncertainty bounds for each $e_i^\text{POD}$ denoted by the triangle markers, where $\sigma_i^\text{POD}$ is the normalized square root of the $i$-th diagonal of the covariance matrix from Equation~\ref{eq:sls-cov}. The results show that the median difference for both sensor configurations is zero, however, the distribution of differences for Configuration 2 has a smaller variance. This again arises from the improved conditioning in the least-squares problem for Configuration 2, resulting in lower uncertainty in the estimated coefficients. From these estimated POD coefficients, we reconstruct the surface pressure field $\hat{\mathbf{p}}$ as in Equation~\ref{eq:pod}. Figure ~\ref{fig:pod_recon_errors} shows the relative L2 errors in the reconstructed pressure field with respect to the reference pressure field, computed by 
\begin{equation}
    e^\text{recon} = \frac{\| \hat{\mathbf{p}} - \mathbf{p}\|}{\| \mathbf{p}\|}.
\end{equation}
We observe that Configuration 1 has a median reconstruction error of 7.5\% in the reconstruction error, but has a large distribution with some errors reaching over 30\%. Configuration 2 has a median of 5.5\% and significantly fewer outliers, with all estimated pressure reconstructions lower than 15\% relative error. 


\begin{figure}[h!]
	\centering
	\includegraphics[width = 6.6in]{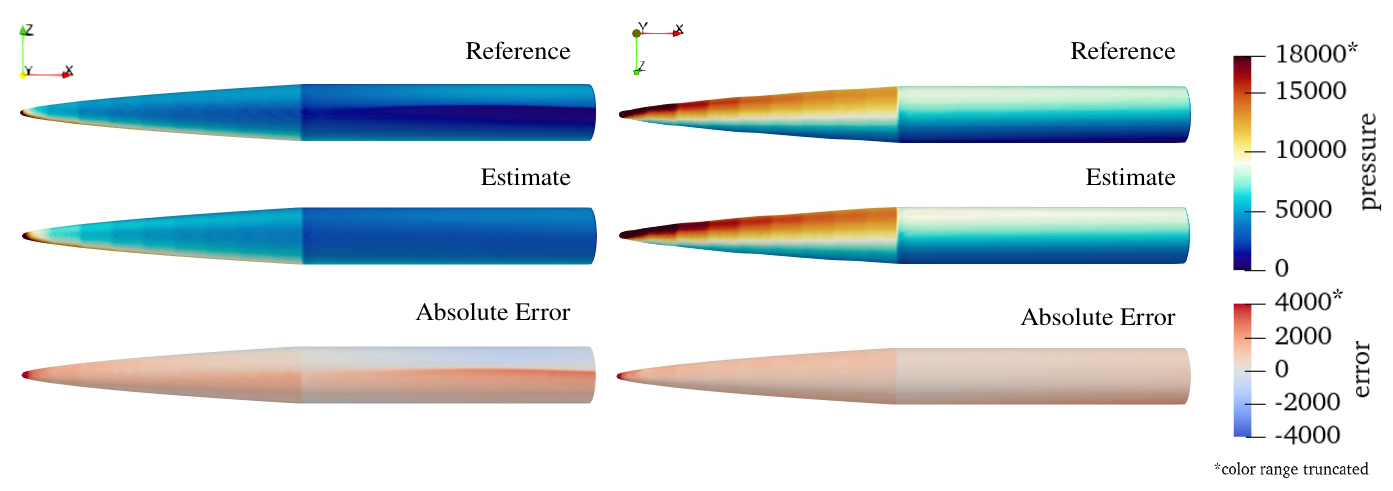}
	\caption{Surface pressure comparison at conditions $M=5, \alpha=6, \beta=5$. (Left) Positive and (right) negative $\mathbf{xz}$-plane view of reference, estimate, and absolute error, from top to bottom.}
	\label{fig:pod_recon_config1}
\end{figure}

\begin{figure}[h!]
	\centering
	\includegraphics[width = 6.6in]{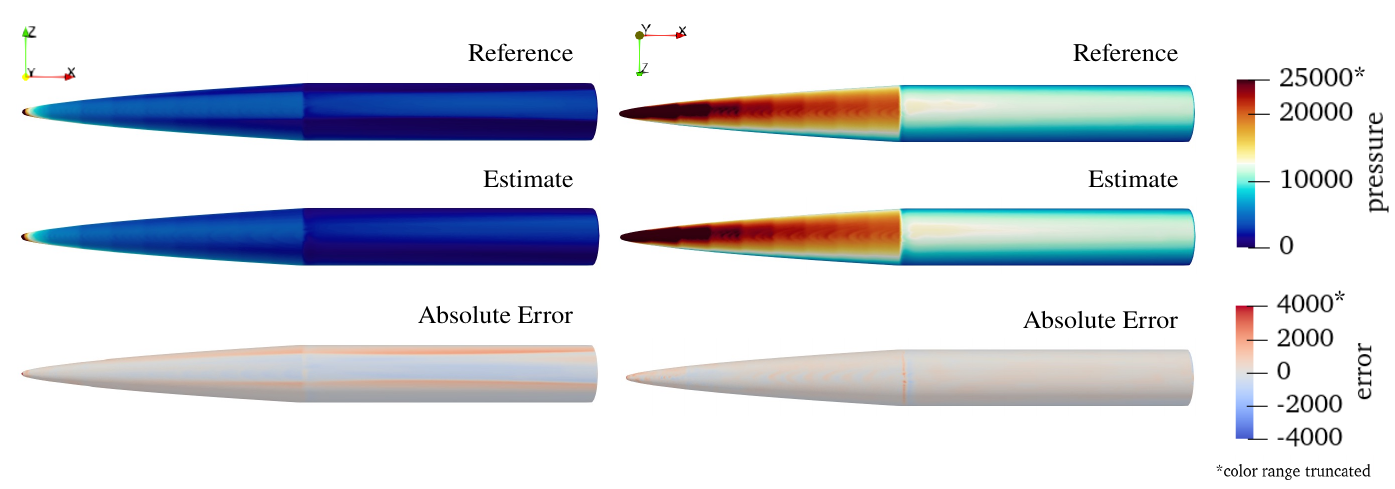}
	\caption{Surface pressure comparison at conditions $M=7, \alpha=2, \beta=8$. (Left) Positive and (right) negative $\mathbf{xz}$-plane view of reference, estimate, and absolute error, from top to bottom.}
	\label{fig:pod_recon_config2}
\end{figure}

To visualize the reconstructed pressure fields, Figure~\ref{fig:pod_recon_config1} and Figure~\ref{fig:pod_recon_config2} each show an example reconstructed surface pressure from the estimated POD coefficients using sensor configuration 1 and 2, respectively. The flight conditions in Figure~\ref{fig:pod_recon_config1} are $M=5$, $\alpha=6$, $\beta=5$, which is in the set $P_1$. The reconstruction error was 15.9\% using a synthetic noise realization for Configuration 1. The flight conditions in Figure~\ref{fig:pod_recon_config2} are $M=7$, $\alpha=2$, $\beta=8$, which is a testing condition not in $P_1$. The reconstruction error was 3.9\% using a synthetic noise realization for Configuration 2. In both cases, the general surface pressure is generally well recovered. We observe that the errors primarily arise from the leeward shock--- the discontinuity is not well captured, which is a result of both noise and the POD surface pressure parameterization, particularly in the testing condition, since the shock location is outside of the range space of the pressure snapshots. 


\begin{figure}[h!]
	\centering
	\includegraphics[width = 5.4in]{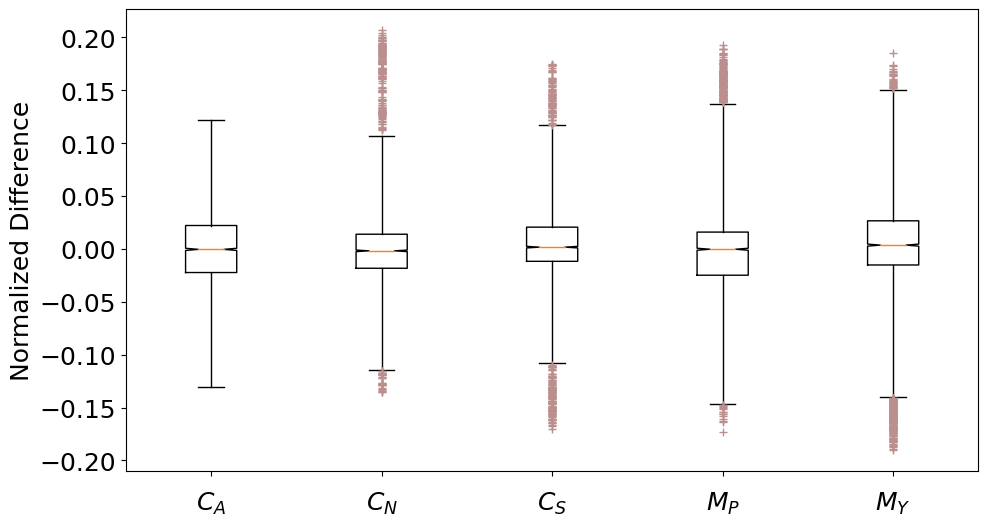}\hspace{1cm}
	\caption{Normalized difference in estimated force and moment coefficients compared to reference values.}
	\label{fig:rel_error_coerr}
\end{figure}

\newpage 
Furthermore, we compute the aerodynamic force and moment coefficients from the estimated POD coefficients. The aerodynamic coefficients, which are the integrated field quantities of interest, can be computed for a pressure field $\mathbf{p}$. We denote the body frame force and moment coefficients by the vector $\mathbf{q}_\text{coeff} = [C_A, C_N, C_Y, M_P, M_Y]^\top$, the elements of which are the coefficients of axial, normal, and side force coefficients, and the pitch and yaw moment coefficients, respectively. These coefficients are computed by
\begin{equation}
    \mathbf{q}_\text{coeff} = \mathbf{G}\mathbf{f} = \mathbf{GC}_\text{map}\mathbf{p} 
\end{equation}
where $\mathbf{G}$ maps the nodal forces to the force and moment coefficients.  For Case 1, since $\mathbf{p} \approx \mathbf{V}_r \mathbf{q} + \bar{\mathbf{p}}$, we can further expand $\mathbf{q}_\text{coeff} =\mathbf{GC}_\text{map}\mathbf{V}_r\mathbf{q} + \mathbf{GC}_\text{map} \bar{\mathbf{p}}$. Note that $\mathbf{GC}_\text{map}\mathbf{V}_r$ and $\mathbf{GC}_\text{map} \bar{\mathbf{p}}$ can be precomputed, allowing rapid computation of the force and moment coefficients in real-time, in the case of guidance and control applications.  Figure~\ref{fig:rel_error_coerr} shows the normalized difference for the estimated force and moment coefficients over the same measurement sample set for Configuration 2, computed by 
\begin{equation}
    e^\text{coeff}_k = \frac{\hat{C}_k - C_k}{\text{max}\{\epsilon_k, |C_k|\}}
\end{equation}
where $C_k$ denotes the $k$-th coefficient in $\mathbf{q}_\text{coeff}$, and $\epsilon_k$ is a tolerance defined as 10\% of the maximum value of $|C_k|$ over all snapshots in $P_1$. This metric is a signed relative error except when the true coefficient is close to zero, where the tolerance $\epsilon_k$ prevents division by zero. We observe that the median normalized difference is close to zero for all coefficients, and the distribution is again a result of the measurement noise. 

\newpage
\subsection{Case 2: Numerical Results}

\begin{figure}[h]
	\centering
	\includegraphics[width = 3.6in]{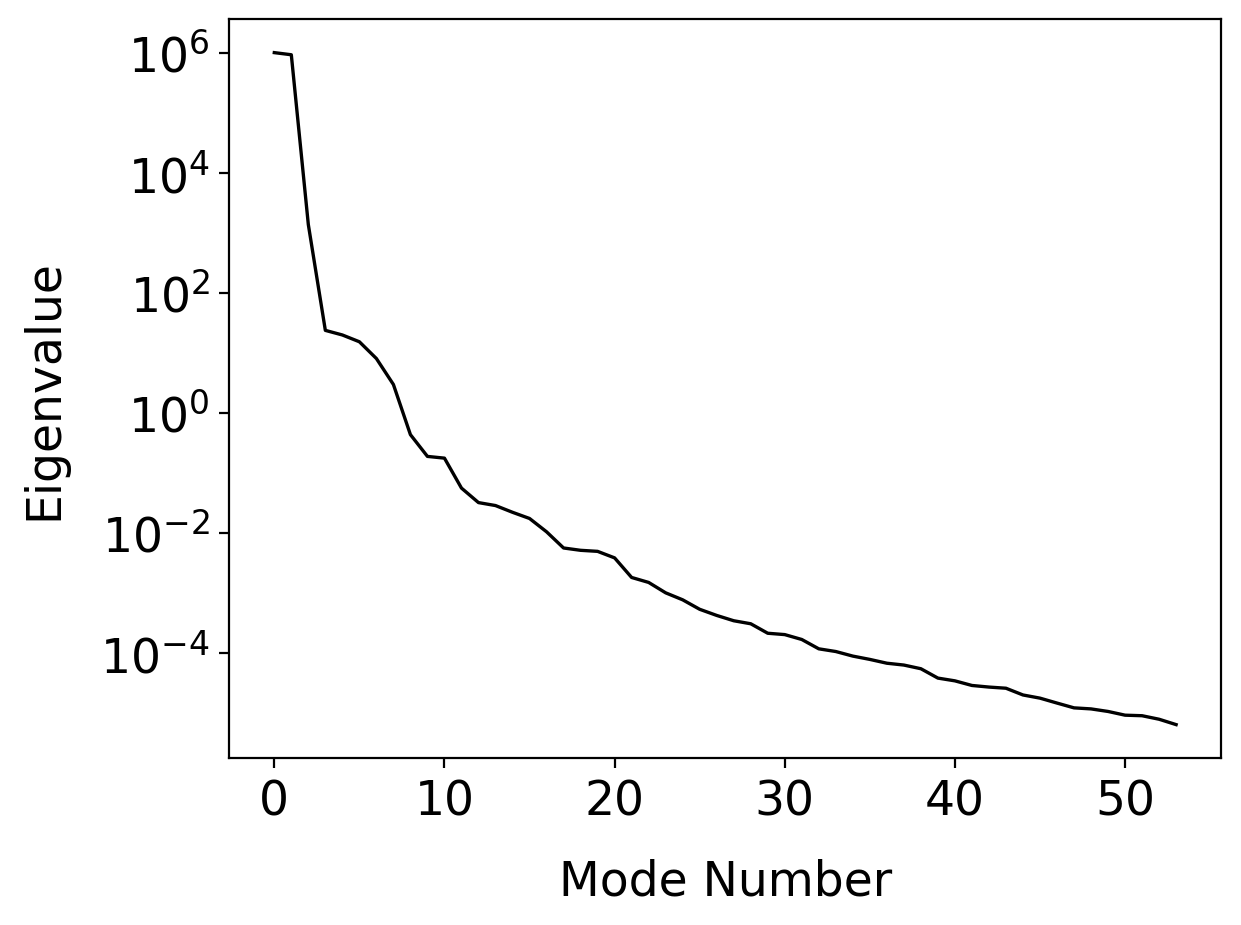}\hspace{1cm}
	\caption{Eigenvalue decay of prior-preconditioned data-misfit Hessian.}
	\label{fig:Heigs}
\end{figure}

\begin{figure}[h!]
	\centering
	\includegraphics[width = 5.5in]{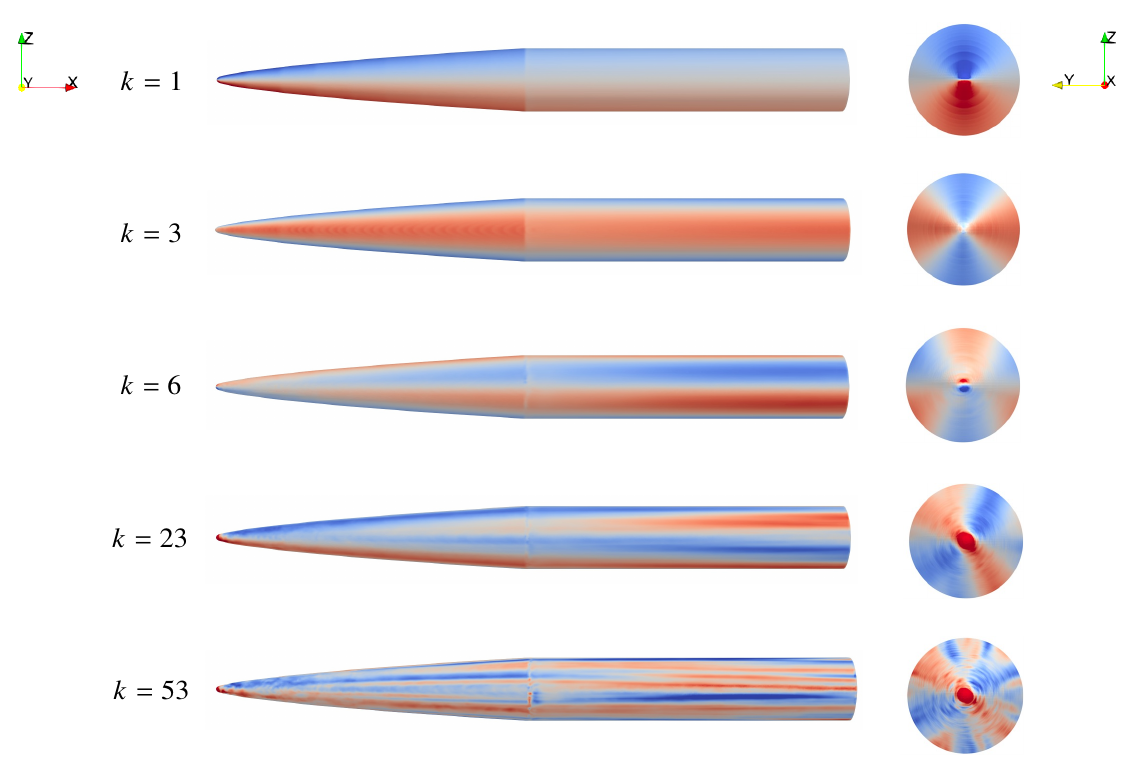}
	\caption{(Left) Side view and (right) front view of the eigenvectors of the prior-preconditioned Hessian for $k=1,3,6,23,53$. }
	\label{fig:Heigvecs}
\end{figure}

For the case where $n_q > n_d$, we seek to estimate the full surface pressure field and quantify the theoretical recoverability from  measurements of the structural strain with uncertainty. Here, we focus on sensor configuration 2 with $n_d = 54$ sensors, since it offers better conditioning and information to assess the limits of the pressure field recovery. We compute the prior mean and covariance using Equation~\ref{eq:prior_cov} for the set of flight conditions given by $P_2$. To analyze the modes which are well-informed by both the prior distribution and the data, we solve the eigenvalue problem with the prior-preconditioned data misfit Hessian, as described in Equation~\ref{eq:eigproblem}. The decay of the non-zero eigenvalues are shown in Figure~\ref{fig:Heigs}. We observe that there is rapid decay in the eigenvalues, with a six order-of-magnitude difference between the first and the eighth eigenvalue. Figure~\ref{fig:Heigvecs} provides a visualization of the  eigenvectors for $k = 1, 3, 6, 23, 53$. We observe that modes corresponding to larger eigenvalues are relatively smooth, and the modes corresponding to smaller eigenvalues contain higher spatial frequency features. These modes represent the features captured by the prior distribution that are also informed by the measurements. However, due to the rapid decay of the eigenvalues, not all of the modes will be recoverable because these eigenvalues will amplify the noise in the measurements. To mitigate this noise amplification, we require tuning of the regularization parameter $\gamma$.

\begin{figure}[h!]
	\centering
	\includegraphics[width = 5.6in]{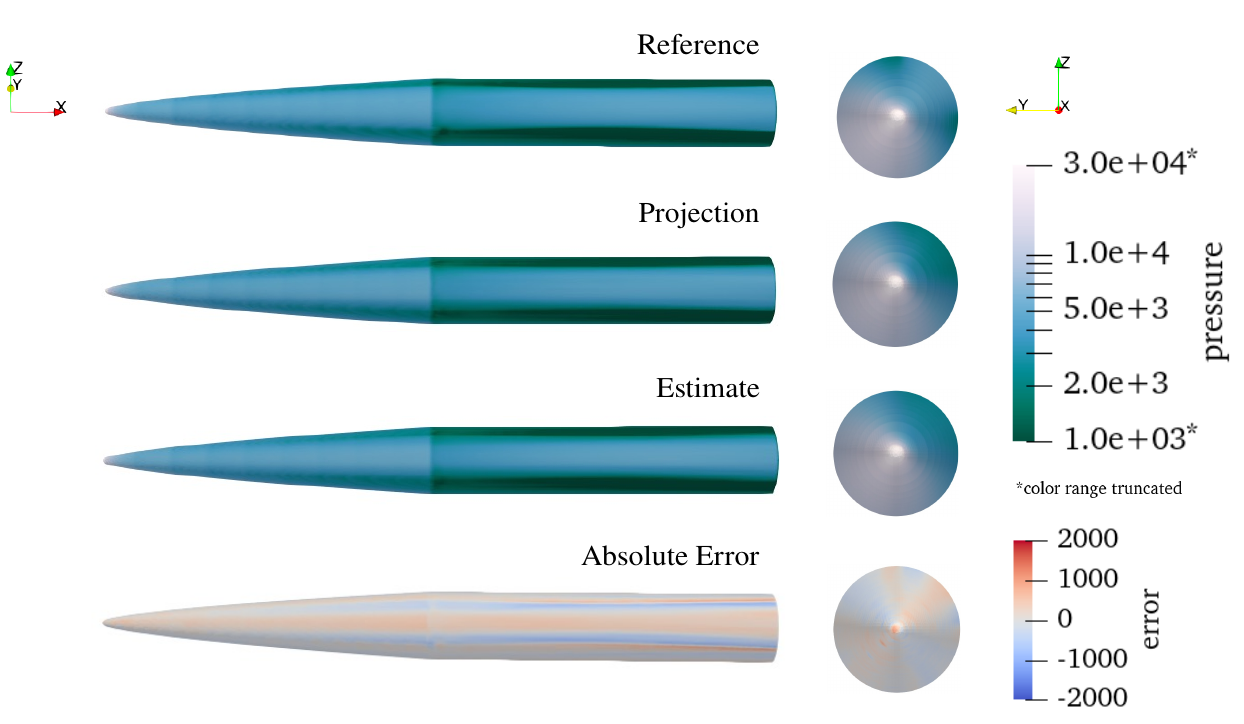}
	\caption{(Left) Leeward view and (right) front view of (top to bottom) the reference pressure, projected reference, zero-noise estimated pressure, and absolute error for $\boldsymbol{\gamma}=10^{-6}$ at flight conditions $M=5, \alpha=6, \beta=6$.}
	\label{fig:pressure_comp_noisefree}
\end{figure}

For a regularization parameter $\gamma$, the contribution of the eigenvectors which correspond to eigenvalues $\lambda_k << \gamma$ will be very small. Consequently, the projection of a centered pressure field $\mathbf{p}$ onto the span of the first $r$ eigenvectors,  choosing $r$ such that $\mathcal{O}(\lambda_r) = \mathcal{O}(\gamma)$ gives an approximate best-case estimate for a given $\gamma$. To demonstrate, we first consider $\gamma = 10^{-6}$ in a noise-free setting to estimate the pressure field at conditions $M=5, \alpha=6, \beta=6$. This is equivalent to retaining all 54 modes in the solution, since all eigenvalues are larger than $\gamma$. Using this parameter, the estimated pressure field is computed for a zero-noise measurement, resulting in a relative error of 2.8\%. Figure~\ref{fig:pressure_comp_noisefree} shows the reference pressure, the projection of the reference onto $r=54$ leading eigenvectors, the estimated pressure, and the absolute error between the reference and the estimate. Here, we observe that the leeward pressure discontinuity is captured by the modes, demonstrating that the strain measurements sufficiently inform this feature of the surface pressure. 

\begin{figure}[h!]
	\centering
	\includegraphics[width = 5.6in]{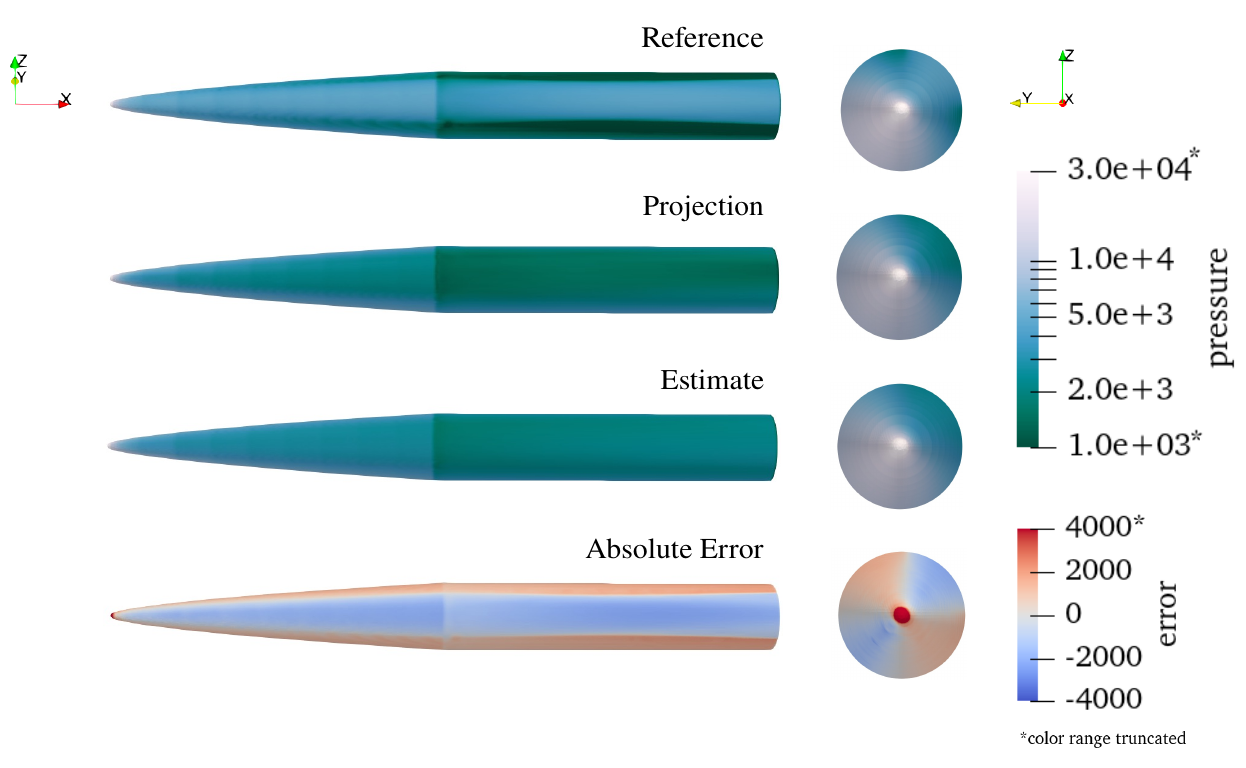}
	\caption{(Left) Leeward view and (right) front view of (top to bottom) the reference pressure, projected reference, estimated pressure, and absolute error for $\boldsymbol{\gamma}=10^{1/2}$ at flight conditions $M=5, \alpha=6, \beta=6$. }
	\label{fig:pressure_comp_reg}
\end{figure}

With the introduction of noise, we must tune $\gamma$ to prevent the noise amplification due to the smaller eigenvalues. For the assumed noise level, the regularization was selected  using a Morozov-like criterion, where the parameter $\gamma$ is selected such that norm of the data misfit between the model strain and the measurement, $\| \mathbf{Z} \mathbf{q} - \mathbf{d} \|$, is approximately equal to the noise level, defined as $\delta = \mathbb{E}[\| \boldsymbol{\eta} \|] = \sigma\sqrt{n_d}$. This resulted in a regularization parameter of $\gamma=10^{1/2}$. We note that a single value for $\gamma$ was selected for estimation over all flight conditions. Considering the same flight conditions as above, $M=5, \alpha=6, \beta=6$, the estimated pressure field was computed for one hundred different noise realizations, resulting in a relative error range of 13.3\% to 27.9\%. Figure~\ref{fig:pressure_comp_reg} shows the reference pressure field, the projection $r=7$, the estimated pressure field for a particular noise realization, which has a relative error of 17.5\%, and the absolute error between the reference and the estimate. We observe that the solution no longer captures the leeward pressure feature because the regularization has eliminated (smoothed) the modes which inform the discontinuity. However, we have sufficiently prevented exploding errors produced by noise amplification. To further assess the estimator performance, we compute the force and moment coefficients for the estimated pressure field. Over the same one hundred synthetic noise realizations for the flight condition above, the relative errors in the coefficient of normal and side force were below 2\%, and the relative errors in the pitch and yaw moment coefficients were below 3.4\%. The largest axial force moment error was larger, reaching 16.8\%, but with a mean of 5.1\% over the one hundred noise realizations. This shows that the integrated quantity relative errors can be quite low in comparison to the full-field relative error. To understand this behavior, returning to Figure~\ref{fig:pressure_comp_reg}, we examine the front view (right column) of the reference and estimated surface pressure for the same flight conditions, as well as the absolute difference. In this view, we can see symmetric structure in the absolute error across the plane with normal $[0, -\cos{\frac{\pi}{4}}, \cos{\frac{\pi}{4}}]$ (direction of the total angle of attack). This results in a similar bending moment in comparison to the reference pressure field, which dominates the strain response. Consequently, the coefficients related to normal and side forces have small relative errors. This behavior arises from the terms in the regularized least-squares problem: we seek a data misfit which is of the same order as the noise level, thus the resulting bending moment is  accurate enough to invoke a similar strain response, however we are unable to correctly recover higher-frequency features in the pressure field due to the regularization for sensor noise.    

\begin{figure}[h!]
	\centering
	\includegraphics[width = 5.6in]{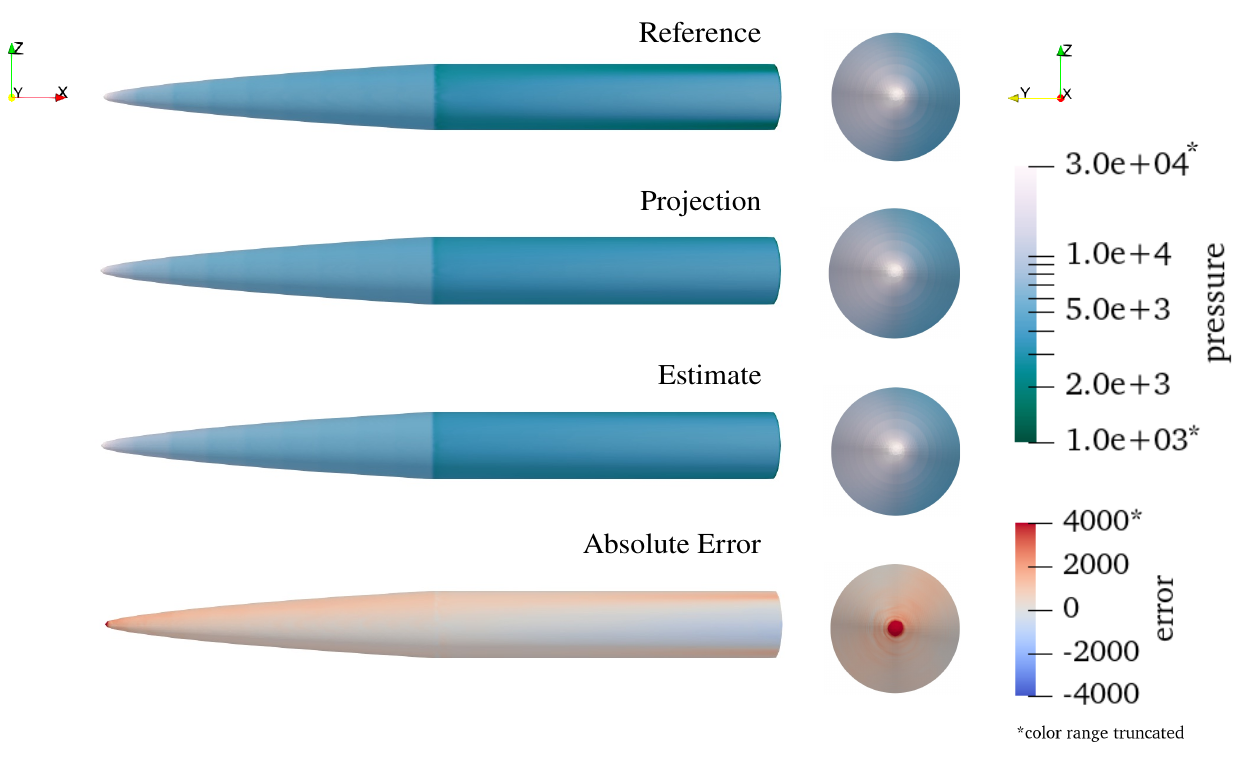}
	\caption{(Left) Side view and (right) front view of (top to bottom) the reference pressure, projected reference, estimated pressure, and absolute error for $\boldsymbol{\gamma}=10^{1/2}$ at testing flight conditions $M=6, \alpha=0, \beta=5$.}
	\label{fig:pressure_comp_reg_tc}
\end{figure}

\newpage
Additionally, we consider flight condition $M=6, \alpha=0, \beta=5$, which was not in the set $P_2$ used to construct the prior covariance. Figure~\ref{fig:pressure_comp_reg_tc} shows the reference pressure field, the projection $r=7$, the estimated pressure field for a particular noise realization, which has a relative error of 7.6\%, and the absolute error between the reference and the estimate. For this flight condition, the leeward discontinuity is less prominent in the reference solution because the total angle of attack is smaller. Again, we observe that the leeward feature is not captured by the estimated pressure field, due to the regularization. Over one hundred noise realizations, the relative error range was 4.0\% up to 12.0\%. This shows consistent estimator performance for flight conditions outside of $P_2$.  It also shows that there is smaller uncertainty in the estimated pressure field for flight conditions at lower angles of attack due to a smoother leeward pressure feature. This highlights the challenge with correctly identifying the discontinuities in a surface pressure field using strain-based inverse maps. 

In summary, the strain-based inverse map using our example sensor configuration captures the structural behavior sufficiently to accurately recover surface pressures over a wide range of flight conditions in a noise-free setting. Therefore, the limitation in pressure recoverability is due to the inherently noisy nature of measurements. We observe that the flexibility of the high-dimensional parameter field in Case 2 still allows the solution to adapt to best match the strain response, resulting in small errors in the performance quantities of interest like the force and moment coefficients. In contrast, the POD parameterization in Case 1 does not enjoy the same advantage, resulting in larger force and moment errors in the presence of sensor noise. In both cases, high spatial frequency features, such as a shock, remain difficult to recover. We note that the pressure recovery presented here can be extended to the vehicle geometry with fins. The recovery of the pressure on the fins, like the fuselage, would be facilitated by the construction of the POD (in Case 1) or the prior covariance (in Case 2). Without the fins, the pressure field still contains complex features such as discontinuities, which is the focus
of our analytical recoverability analysis. The surface pressure on the fins would not contain these features, though they could potentially cause a secondary shock footprint on the fuselage near the aft end. Our analysis in this work already analyzes the recoverability of these types of challenging characteristics in the presence of sensor noise, and can be applied in the same fashion for finned vehicle geometries.

\section{Conclusion}\label{sec:conclusion}

In this work, we have formulated a real-time estimator for aerodynamic loads via a strain-based strategy. Pre-computation of the high-dimensional system matrix multiplications in the resultant estimator, which embed the physics in the solution, enables real-time tractability of the approach. Given sensor noise statistics, the approach also enables uncertainty quantification for the estimated quantities of interest. In the case where the number of unknown pressure quantities is smaller than the number of measurements (i.e., $n_q \leq n_d$), we directly compute the variance of the estimated parameters, demonstrating the impacts of the conditioning of the parameter-to-observable map. In the case where the number of pressure quantities exceeds the number of measurements (i.e., $n_q > n_d)$, we analyze the prior-preconditioned Hessian to expose the theoretical recoverability of pressure fields under measurement uncertainty through the prior distribution and the data. Although two different sensor configurations were considered in this paper, the question of optimal sensor placement for strain-based sensing strategies will be addressed in future work. Future work will also consider aerothermal heating effects, which impact the strain response of the vehicle. The developments in this paper, as well as the planned studies, are critical steps towards a deployable strain-based aerodynamic measurement strategy for controls, ground and flight testing in hypersonics. 

\section*{Acknowledgments}
This work was supported by AFOSR grant FA9550-21-1-0089 under the NASA University Leadership Initiative (ULI), the DOE ASCR grant DE-SC0021239, and the DOD MURI grant FA9550-21-1-0084. This material is based on work supported by the National Science Foundation Graduate Research Fellowship under Grant No. DGE 2137420. The authors would like to thank the members of the FAST ULI team for the engaging discussions and support in this work. 

\newpage
\bibliographystyle{ieeetr}
\bibliography{ref}

\end{document}